\documentclass[11pt]{amsart}

\usepackage{xspace}

\newcommand{\1}{{1\!\!1}}

\usepackage{amssymb}

\makeatletter

\newdimen\@Leftmargin   \@Leftmargin=3.0cm
\newdimen\@Rightmargin  \@Rightmargin=3.0cm
\newdimen\@Topmargin    \@Topmargin=2.0cm
\newdimen\@Bottommargin \@Bottommargin=2.0cm
\def\InitLayout{
  \setlength{\textwidth}{\paperwidth}
  \addtolength{\textwidth}{-\@Leftmargin}
  \addtolength{\textwidth}{-\@Rightmargin}
  \setlength{\textheight}{\paperheight}
  \addtolength{\textheight}{-\@Topmargin}
  \addtolength{\textheight}{-\@Bottommargin}
  \addtolength{\textheight}{-\headheight}
  \addtolength{\textheight}{-\headsep}
  \addtolength{\textheight}{-\footskip}
  \setlength{\oddsidemargin}{\@Leftmargin}
  \addtolength{\oddsidemargin}{-1in}
  \setlength{\evensidemargin}{\@Rightmargin}
  \addtolength{\evensidemargin}{-1in}
  \setlength{\topmargin}{\@Topmargin}
  \addtolength{\topmargin}{-1in}
}
\InitLayout

\let\c@equation=\c@subsection

\def\@sect#1#2#3#4#5#6[#7]#8{%
 \edef\@toclevel{\ifnum#2=\@m 0\else\number#2\fi}%
 \@tempskipa #5\relax
 \ifnum #2>\c@secnumdepth
  \ifdim \@tempskipa>\z@ 
   \begingroup #6\relax
   \@hangfrom{\hskip #3\relax}{\interlinepenalty\@M #8\par}%
   \endgroup
  \else
   \def\@svsechd{#6\hskip #3
   \@ifnotempty{#8}{\ignorespaces#8\unskip\@addpunct.}%
   }%
  \fi
 \else
  \@xp\let\@xp\@secnumber\csname the#1\endcsname%
  \refstepcounter{#1}%
  \edef\@svsec{\ifnum#2<\@m\@ifundefined{#1name}{}{\csname #1name\endcsname}\fi
  \@nx\textup{\@nx\mdseries\csname the#1\endcsname.\enspace}}
  \ifdim \@tempskipa>\z@ 
   \begingroup #6\relax
   \@hangfrom{\hskip #3\relax\@svsec}{\interlinepenalty\@M #8\par}%
   \endgroup
   \ifnum#2>\@m \else \@tocwrite{#1}{#8}\fi
  \else
   \def\@svsechd{#6\hskip #3\@svsec\enspace
   \@ifnotempty{#8}{\ignorespaces#8\unskip\@addpunct.}%
   \ifnum#2>\@m \else \@tocwrite{#1}{#8}\fi
   }%
  \fi
 \fi
 \global\@nobreaktrue
 \@xsect{#5}}

\makeatother

\def\to{\mathchoice
{\longrightarrow}
{\rightarrow}
{\rightarrow}
{\rightarrow}}

\def\mapsto{\DOTSB\mapstochar\to}

\def\hookrightarrow{\mathchoice
{\DOTSB\lhook\joinrel\relbar\joinrel\rightarrow}
{\DOTSB\lhook\joinrel\rightarrow}
{\DOTSB\lhook\joinrel\rightarrow}
{\DOTSB\lhook\joinrel\rightarrow}}

\newtheorem{theorem}[subsection]{Theorem}
\newtheorem*{theorem*}{Theorem}
\newtheorem{proposition}[subsection]{Proposition}
\newtheorem{lemma}[subsection]{Lemma}
\newtheorem{corollary}[subsection]{Corollary}

\theoremstyle{definition}

\newtheorem{definition}[subsection]{Definition}
\newtheorem{example}[subsection]{Example}

\newcommand{\Z}{\mathbb{Z}}
\newcommand{\C}{\mathbb{C}}
\newcommand{\N}{\mathbb{N}}

\newcommand{\Q}{\mathbb{Q}}

\newcommand{\eps}{\varepsilon}
\renewcommand{\o}{\otimes}
\newcommand{\ohat}{\Hat{\otimes}}

\DeclareMathOperator{\PSL}{PSL}
\DeclareMathOperator{\gl}{\mathfrak{gl}}

\DeclareMathOperator{\Ind}{Ind}

\DeclareMathOperator{\End}{End}

\newcommand{\p}{\partial}

\newcommand{\CM}{\mathcal{M}}

\DeclareMathOperator{\gr}{gr}
\renewcommand{\*}{\cdot}

\newcommand{\<}{\langle}
\renewcommand{\>}{\rangle}

\newcommand{\Mbar}{\overline{\mathcal{M}}}
\renewcommand{\]}{{]\!]}}
\renewcommand{\[}{{[\![}}
\DeclareMathOperator{\ch}{ch}

\DeclareMathOperator{\Tr}{Tr}
\DeclareMathOperator{\Exp}{Exp}
\DeclareMathOperator{\Log}{Log}
\newcommand{\CP}{\mathbb{CP}}
\newcommand{\bull}{\bullet}
\renewcommand{\SS}{\mathbb{S}}

\DeclareMathOperator{\spec}{spec}
\newcommand{\point}{\mathbb{A}^0}
\renewcommand{\P}{\mathsf{P}}
\newcommand{\Po}{\overset{\circ}{\mathsf{P}}}

\newcommand{\TT}{{\mathbb{T}}}

\DeclareMathOperator{\GS}{\mathsf{E}}
\DeclareMathOperator{\Serre}{\mathsf{e}}
\newcommand{\VV}{\mathbb{V}}

\newcommand{\Proj}{{\mathsf{Proj}}}
\newcommand{\Var}{{\mathsf{Var}}}
\newcommand{\Mot}{{\textsf{Mot}}}

\DeclareMathOperator{\im}{Im}

\newcommand{\ring}{rring\xspace}
\newcommand{\rings}{rrings\xspace}
\DeclareMathOperator{\FM}{\mathsf{FM}}
\newcommand{\FF}{\mathsf{F}}
\newcommand{\T}{\mathsf{T}}
\newcommand{\G}{{\mathbb{G}}}
\newcommand{\zbar}{{\bar{z}}}
\renewcommand{\hom}[2][{}]{{}^{#1}\[\SS,#2\]}
\newcommand{\ZZ}{{\mathsf{Z}}}

\DeclareMathOperator{\VERT}{Vert}

\newcommand{\V}{\mathcal{V}}
\newcommand{\w}{\mathcal{W}}
\newcommand{\x}{\mathcal{X}}
\newcommand{\y}{\mathcal{Y}}
\newcommand{\LL}{\mathsf{L}}

\newcommand{\CR}{\mathcal{R}}

\newcommand{\CS}{\mathcal{S}}

\def\({(\!(}
\def\){)\!)}

\begin{document}

\title{Mixed Hodge structures of configuration spaces}

\author{E. Getzler}

\address{Department of Mathematics, MIT, Cambridge MA 02139 USA}

\email{getzler@math.mit.edu}

\maketitle

Let $X$ be a smooth projective variety over $\C$. The configuration space
$\FF(X,n)$ of $X$ is the complement of the diagonals in $\T(X,n)=X^n$:
$$
\FF(X,n) = \{ (z_1,\dots,z_n) \in X^n \mid \text{$z_i\ne z_j$ for $i\ne j$}
\} .
$$
The symmetric group $\SS_n$ acts freely on $\FF(X,n)$; in this paper, we
study the induced action of the symmetric group $\SS_n$ on
$H^{p,q}(\FF(X,n))$. As an application of our results, we calculate the
$\SS_n$-equivariant Hodge polynomial of the Fulton-MacPherson
compactification of $\FF(X,n)$.

In a sequel to this paper, we extend our results to the relative setting:
this is technically more difficult, requiring Saito's theory of mixed Hodge
modules. As an application of the relative theory, we will calculate the
$\SS_n$-equivariant Hodge polynomials of the projective varieties
$\Mbar_{1,n}$. (We have calculated the $\SS_n$-equivariant Hodge
polynomials of the projective varieties $\Mbar_{0,n}$ in \cite{gravity}.)

Let $X$ be a complex quasi-projective variety. The Serre polynomial
$\Serre(X)$ of $X$ is the polynomial in variables $u$ and $v$, satisfying
(and indeed characterized by) the following axioms:
\begin{enumerate}
\item if $X$ is projective and smooth, $\Serre(X)$ is the Hodge polynomial%
\footnote{Here, we have modified the usual conventions, replacing $u$ and
$v$ by $-u$ and $-v$. This will lead to cleaner formulas later.}
$$
\Serre(X) = \sum_{p,q=0}^\infty (-u)^p (-v)^q \dim H^{p,q}(X,\C) ;
$$
\item if $Z$ is a closed subvariety of $X$, then
$\Serre(X)=\Serre(X\setminus Z)+\Serre(Z)$.
\end{enumerate}

A formula for the Serre polynomial follows from mixed Hodge theory (Deligne
\cite{Deligne}):
\begin{equation} \label{Serre}
\Serre(X) = \sum_{p,q=0}^\infty u^p v^q
\chi\bigl( H^\bull_c(X,\C)^{p,q} \bigr) , 
\end{equation}
where if $(V,F,W)$ is a mixed Hodge structure over $\C$,
$$
V^{p,q}=F^p\gr^W_{p+q}V\cap\bar{F}^q\gr^W_{p+q}V ,
$$
and $\chi(V)$ denotes the Euler characteristic of the finite-dimensional
graded vector space $V$. This formula was introduced by Danilov and
Khovanski\v\i, although they do not give it a name: earlier, Serre had
observed that the Weil conjectures, together with resolution of
singularities, implied the existence of such a polynomial.

The Serre polynomial represents the class of $H^\bull_c(X,\C)$ in the
Grothendieck group of mixed Hodge structures over $\C$. It is a
``character'' on varieties, since it satisfies the K\"unneth formula:
$$
\Serre(X\times Y) = \Serre(X) \Serre(Y) .
$$
We borrow from Manin \cite{Manin} the notation $\LL$ for the Serre
polynomial $uv$ of $\C(-1)=H^2(\CP^1,\C)$.

Denote by $S^nX$ the $n$th symmetric power of the variety $X$, that is, the
quotient of $X^n$ by the symmetric group $\SS^n$. It follows from the
formula of Macdonald \cite{Macdonald-symmetric} that the generating
function $\sigma_t(X)$ for the Serre polynomials $\Serre(S^nX)$ has the
formula
$$
\sigma_t(X) = \sum_{n=0}^\infty t^n \Serre(S^nX) = \prod_{p,q}
(1-tu^pv^q)^{-\chi(H^\bull_c(X,\C)^{p,q})} .
$$
(This is the Hodge analogue of the zeta function of a variety over a finite
field.)

If a finite group $G$ acts on $X$, define the equivariant Serre polynomial
by the formula
$$
\Serre_g(X) = \sum_{p,q=0}^\infty u^pv^q \sum_i (-1)^i
\Tr(g|H^i_c(X,\C)^{p,q}) .
$$

Let $\FF(\C,n)$ be the configuration space of $n$ ordered points in
$\C$. By a formula of Lehrer and Solomon \cite{LS} (see
\eqref{Lehrer-Solomon}), we see that if $\sigma\in\SS_n$ has $n_j$ cycles
of length $j$, then
$$
\Serre_\sigma(\FF(\C,n)) = \prod_{j=1}^\infty \alpha_j \bigl( \alpha_j -
jt^j \bigl) \dots \bigl( \alpha_j - (n_j-1)jt^j \bigl) ,
$$
where $\alpha_j=\sum_{d|j}\mu(j/d)\LL^d$. (For another proof of this
formula, see \cite{gravity}.) In this paper, we prove the following
generalization of this formula.
\begin{theorem*}
Let $X$ be a quasi-projective variety, and let
$\alpha_j(X)=\sum_{d|j}\mu(j/d)\Serre(X;u^d,v^d)$. If $\sigma\in\SS_n$ is a
permutation with $n_j$ cycles of length $j$, then
$$
\Serre_\sigma(\FF(X,n)) = \prod_{j=1}^\infty \alpha_j(X) \bigl( \alpha_j(X)
- j \bigl) \dots \bigl( \alpha_j(X) - (n_j-1)j \bigl) .
$$
\end{theorem*}

Given a $\SS_n$-module $V$, consider the local system
$\VV=\FF(X,n)\times_{\SS_n}V$ over $\FF(X,n)/\SS_n$. This local system has
a Serre polynomial
$$
\Serre(X,V) = \sum_{p,q=0}^\infty u^p v^q \chi\bigl( H^\bull_c(X,\VV)^{p,q}
\bigr) ,
$$
related to the $\SS_n$-equivariant Serre polynomial of $\FF(X,n)$ by the
formula
$$
\Serre(X,V) = \frac{1}{n!} \sum_{\sigma\in\SS_n} \Tr(\sigma|V)
\Serre_\sigma(\FF(X,n)) .
$$
A special case of this is the trivial representation $\1$, for which
$\Serre(\FF(X,n),\1)=\Serre(\FF(X,n)/\SS_n)$: we will prove that
$$
\sum_{n=0}^\infty x^n \Serre(\FF(X,n)/\SS_n) =
\frac{\sigma_t(X)}{\sigma_{t^2}(X)} = \prod_{p,q} \Bigl(
\frac{1-t^2u^pv^q}{1-tu^pv^q} \Bigr)^{\chi(H_c^\bull(X,\C)^{p,q})} .
$$
The analogue of this formula for $X=\spec(\Z)$ is the Dirichlet series
$$
\frac{\zeta(s)}{\zeta(2s)} = \sum_{n=1}^\infty \frac{\mu(n)^2}{n^s} ,
$$
where $\mu(n)^2$ is the arithmetic function which is $1$ if $n$ is
square-free and $0$ otherwise: the analogy with configuration spaces is
clear.

Although the above formulas involve de Rham cohomology, they may be proved
with no greater difficulty for motivic cohomology, using the recent results
of Gillet and Soul\'e \cite{GS}. Let $\Mot$ be the category of (pure
effective rational) Chow motives, with Grothendieck group $K_0(\Mot)$.
Given a quasi-projective variety $X$ over $\C$, let $\Serre(X)\in
K_0(\Mot)$ be the virtual motive associated to $X$ by Gillet and
Soul\'e. Then $K_0(\Mot)$ is a $\lambda$-ring, with $\sigma$-operations
satisfying
$$
\sigma_n(\Serre(X)) = \Serre(X^n/\SS_n) ,
$$
and $\lambda$-operations satisfying
$$
\lambda_n(\Serre(X)) = \Serre(\FF(X,n),\eps) ,
$$
where $\eps$ is the sign representation of $\SS_n$. For each irreducible
representation $V_\lambda$ of $\SS_n$, we will construct a natural
transformation $\Phi_\lambda$ of the category of $\lambda$-rings such that
$$
\Phi_\lambda(\Serre(X)) = \Serre(\FF(X,n),V_\lambda) ,
$$
where $V_\lambda$ is the irreducible $\SS_n$-module associated to the
partition $\lambda$. Here is a table of $\Phi_\lambda$ for small $\lambda$:
$$\begin{tabular}{l||l}
$\Phi_1$	& $\sigma_1$ \\ \hline
$\Phi_2$	& $\sigma_2 - \sigma_1$ \\
$\Phi_{1^2}$	& $\sigma_{1^2}$ \\ \hline
$\Phi_3$	& $\sigma_3 - \sigma_2 - \sigma_{1^2}$ \\
$\Phi_{21}$	& $\sigma_{21} - \sigma_2 - \sigma_{1^2} + \sigma_1$ \\
$\Phi_{1^3}$	& $\sigma_{1^3}$ \\ \hline
$\Phi_4$	& $\sigma_4 - \sigma_3 - \sigma_{21} + \sigma_{1^2}$ \\
$\Phi_{31}$	& $\sigma_{31} - \sigma_3 - 2\sigma_{21} - \sigma_{1^3} +
			2\sigma_2 + \sigma_{1^2} - \sigma_1$ \\
$\Phi_{2^2}$	& $\sigma_{2^2} - \sigma_3 - \sigma_{21} + \sigma_2 +
			2\sigma_{1^2}$ \\
$\Phi_{21^2}$	& $\sigma_{21^2} - \sigma_{21} - \sigma_{1^3} + \sigma_2
			+ \sigma_{1^2} - \sigma_1$ \\
$\Phi_{1^4}$	& $\sigma_{1^4}$ \\
\end{tabular}$$
(Here, $\sigma_\lambda$ denotes the operation on $\lambda$-rings associated
to the Schur polynomial $s_\lambda$.) For the generating function of the
operations $\Phi_\lambda$, see \eqref{Phi}.

The organization of this paper is as follows. In Section 1, we recall some
of the theory of symmetric functions in an infinite number of variables,
and its relationship to the theory of $\lambda$-rings. In Section 2, we
introduce the completion of a $\lambda$-ring: this is used in the proof of
our main result, where we work with generating (symmetric) functions, which
lie in such a completion. We also introduce the $\Phi$-operations, which
enter into the statement of our formula for $\Serre(\FF(X,n),V_\lambda)$.

In Section 3, we introduce a class of categories, Karoubian \rings (sic),
which have many of the properties of the category of modules over a
ring. In particular, we prove the Peter-Weyl Theorem: any representation of
a finite group in a Karoubian \ring over a field of characteristic zero is
completely reducible.  Similar results have been obtained by del Ba\~no
Rollin \cite{Rollin}. In Section 4, we apply this result to construct a
$\lambda$-ring structure on the Grothendieck group of a Karoubian \ring
(over a field of characteristic zero). Section 5 contains our main result,
Theorem \ref{MAIN}, which is a theorem about Serre functors; this is the
name we give to a sequence $\{\x\mapsto\GS^n(\x)\mid n\ge0\}$ of functors
from the category of quasi-projective varieties to a Karoubian \ring $\CR$,
satisfying appropriate analogues of the Meyer-Vietoris and the K\"unneth
theorems. Examples of Serre functors are the de Rham cohomology and the
cohomology of the weight complex of Gillet and Soul\'e \cite{GS}.

In Section 6, we give a simple application of these results, to the
calculation of the $\SS_n$-equivariant Hodge polynomial of the
Fulton-Macpherson compactification of the configuration space $\FF(X,n)$.

\subsection*{Acknowledgments} The author wishes to thanks the Department of
Mathematics at the Universit\'e de Paris-VII for its hospitality during the
inception of this paper. He is partially supported by a research grant of
the NSF, a fellowship of the A.P. Sloan Foundation, and the
Max-Planck-Institut f\"ur Mathematik in Bonn.


\section{Symmetric functions and $\lambda$-rings}

\subsection{Symmetric functions} In this section, we recall some results on
symmetric functions and representations of $\SS_n$ which we need later. For
the proofs of these results, we refer to Macdonald \cite{Macdonald}.

The ring of symmetric functions is the inverse limit
$$
\Lambda = \varprojlim \Z[x_1,\dots,x_k]^{\SS_k} .
$$
It is is a polynomial ring in the complete symmetric functions
$$
h_n = \sum_{i_1\le\dots\le i_n} x_{i_1}\dots x_{i_n} .
$$
The power sums (also known as Newton polynomials)
$$
p_n = \sum_i x_i^n
$$
form a set of generators of the polynomial ring $\Lambda_\Q=\Lambda\o\Q$.
This is shown by means of the elementary formula
\begin{equation} \label{P-H}
P_t = t \frac{d}{dt} \log H_t ,
\end{equation}
where
$$
H_t = \sum_{n=0}^\infty t^n h_n = \prod_i (1-tx_i)^{-1}
\quad\text{and}\quad
P_t = \sum_{n=0}^\infty t^n p_n = \sum_i (1-tx_i)^{-1} .
$$
Written out explicitly, we obtain Newton's formula relating the two sets of
generators:
$$
nh_n = p_n + h_1p_{n-1} + \dots + h_{n-1}p_1 .
$$
We may also invert \eqref{P-H}, obtaining the formula
\begin{equation} \label{H-P}
H_t = \exp \Bigl( \sum_{n=1}^\infty \frac{t^np_n}{n} \Bigr) .
\end{equation}

A partition $\lambda$ is a decreasing sequence
$(\lambda_1\ge\dots\ge\lambda_\ell)$ of positive integers; we write
$\lambda\vdash n$, where $n=\lambda_1+\dots+\lambda_\ell$, and denote the
length of $\lambda$ by $\ell(\lambda)$. Identifying $\Lambda$ with the ring
of characters of the Lie algebra $\gl_\infty = \varprojlim \gl_k$, we see
that partitions correspond to dominant weights, and thus $\Lambda$ has a
basis of consisting of the characters of the irreducible representations of
$\gl_\infty$. These characters, given by the Weyl character formula
$$
s_\lambda = \lim_{k\to\infty}
\frac{\det(x_i^{\lambda_j+k-j})_{1\le i,j\le k}}
{\det(x_i^{k-j})_{1\le i,j\le k}} ,
$$
are known as the Schur functions. In terms of the polynomial generators
$h_n$, they are given by the Jacobi-Trudy formula
$s_\lambda=\det\bigl(h_{\lambda_i-i+j}\bigr)_{1\le i,j\le\ell(\lambda)}$.

There is a non-degenerate integral bilinear form on $\Lambda$, denoted
$\<f,g\>$, for which the Schur functions $s_\lambda$ form an orthonormal
basis. (This is sometimes called the Hall inner product.) The adjoint of
multiplication by $f\in\Lambda$ with respect to this inner product is
denoted $D(f)$. Written in terms of the power sums $p_n$, the operator
$D(f)$ has the formula $D(p_n)=n\p/\p p_n$, while the inner product
$\<f,g\>$ has the formula
$$
\< f , g \> = D(f)g \Big|_{p_n=0,n\ge1} .
$$

\subsection{Pre-$\lambda$-rings} A pre-$\lambda$-ring is a commutative ring
$R$, together with a morphism of commutative rings $\sigma_t:R\to R\[t\]$
such that $\sigma_t(a)=1+ta+O(t^2)$. Expanding $\sigma_t$ in a power series
$$
\sigma_t(a) = \sum_{n=0}^\infty t^n \sigma_n(a) ,
$$
we obtain endomorphisms $\sigma_n$ of $R$ such that $\sigma_0(a)=1$,
$\sigma_1(a)=a$, and
$$
\sigma_n(a+b) = \sum_{k=0}^n \sigma_{n-k}(a) \sigma_k(b) .
$$
There are also operations $\lambda_k(a)=(-1)^k\sigma_k(-a)$, with
generating function
\begin{equation} \label{invert}
\lambda_t(a) = \sum_{n=0}^\infty t^n \lambda_n(a) = \sigma_{-t}(a)^{-1} .
\end{equation}
The $\lambda$-operations are polynomials in the $\sigma$-operations with
integral coefficients, and vice versa. In this paper, we take the
$\sigma$-operations to be more fundamental; nevertheless, following custom,
the structure they define is called a pre-$\lambda$-ring.

Given a pre-$\lambda$-ring $R$, there is a bilinear map $\Lambda\o R\to R$,
which we denote $f\circ a$, defined by the formula
$$
(h_{n_1}\dots h_{n_k})\circ a = \sigma_{n_1}(a)\dots\sigma_{n_k}(a) .
$$
The image of the power sum $p_n$ under this map is the operation on $R$
known as the Adams operation $\psi_n$.  We denote the operation
corresponding to the Schur function $s_\lambda$ by $\sigma_\lambda$. Note
that \eqref{H-P} implies the relation
$$
\sigma_t(a) = \exp \Bigl( \sum_{n=1}^\infty \frac{t^n\psi_n(a)}{n} \Bigr) .
$$

The following formula (I.4.2 of \cite{Macdonald}) is known as Cauchy's
formula:
\begin{equation} \label{Cauchy}
H_t(...,x_iy_j,...) = \prod_{i,j} (1-tx_iy_j)^{-1}
= \sum_{\lambda\vdash n} s_\lambda(x) \o s_\lambda(y)
= \exp\Bigl( \sum_{k=1}^\infty \frac{p_k(x) \o p_k(y)}{k} \Bigr) .
\end{equation}
From it, the following result is immediate.
\begin{proposition}
If $R$ and $S$ are pre-$\lambda$-rings, their tensor product $R\o S$ is a
pre-$\lambda$-ring, with $\sigma$-operations
$$
\sigma_n(a\o b) = \sum_{\lambda\vdash n} \sigma_\lambda(a) \o
\sigma_\lambda(b) ,
$$
and Adams operations $\psi_n(a\o b) = \psi_n(a) \o \psi_n(b)$.
\end{proposition}

For example, $\sigma_2(a\o b) = \sigma_2(a)\o\sigma_2(b) +
\lambda_2(a)\o\lambda_2(b)$.

\subsection{$\lambda$-rings} The polynomial ring $\Z[x]$ is a
pre-$\lambda$-ring, with $\sigma$-operations characterized by the formula
$\sigma_n(x^i)=x^{ni}$. Taking tensor powers of this pre-$\lambda$-ring
with itself, we see that the polynomial ring $\Z[x_1,\dots,x_k]$ is a
pre-$\lambda$-ring. The $\lambda$-operations on this ring are equivariant
with respect to the permutation action of the symmetric group $\SS_k$ on
the generators, hence the ring of symmetric functions
$\Z[x_1,\dots,x_k]^{\SS_k}$ is a pre-$\lambda$-ring. Taking the limit
$k\to\infty$, we obtain a pre-$\lambda$-ring structure on $\Lambda$.

\begin{definition}
A $\lambda$-ring is pre-$\lambda$-ring such that if $f,g\in\Lambda$ and
$x\in R$,
\begin{equation} \label{lambda-ring}
f\circ(g\circ x)=(f\circ g)\circ x .
\end{equation}
\end{definition}

By definition, the pre-$\lambda$-ring $\Lambda$ is a $\lambda$-ring; in
particular, the operation $f\circ g$, called plethysm, is associative.

The following result (see Knutson, \cite{Knutson}) is the chief result in
the theory of $\lambda$-rings.
\begin{theorem} \label{universal}
$\Lambda$ is the universal $\lambda$-ring on a single generator $h_1$.
\end{theorem}

This theorem makes it straighforward to verify identities in
$\lambda$-rings: it suffices to verify them in $\Lambda$. As an
application, we have the following corollary.
\begin{corollary}
The tensor product of two $\lambda$-rings is a $\lambda$-ring.
\end{corollary}
\begin{proof}
We need only verify this for $R=\Lambda$. A torsion-free pre-$\lambda$-ring
whose Adams operations are ring homomorphisms which satisfy
$\psi_m(\psi_n(a))=\psi_{mn}(a)$ is a $\lambda$-ring.  It is easy to verify
these conditions for $\Lambda\o\Lambda$, since $\psi_n(a\o b) = \psi_n(a)
\o \psi_n(b)$.
\end{proof}

In the definition of a $\lambda$-ring, it is usual to adjoin the axiom
$$
\sigma_n(xy) = \sum_{\lambda\vdash n} \sigma_\lambda(a) \o
\sigma_\lambda(y) .
$$
However, this formula follows from our definition of a $\lambda$-ring: by
universality, it suffices to check it for $R=\Lambda\o\Lambda$, $x=h_1\o1$
and $y=1\o h_1$, for which it is evident.

\section{Complete $\lambda$-rings}

A filtered $\lambda$-ring $R$ is a $\lambda$-ring with decreasing
filtration
$$
R = F^0R \supset F^1R \supset \dots ,
$$
such that
\begin{enumerate}
\item $\bigcap_k F^kR = 0$ (the filtration is discrete);
\item $F^mR\*F^nR\subset F^{m+n}R$ (the filtration is compatible with the
product);
\item $\sigma_m(F^nR)\subset F^{mn}R$ (the filtration is compatible with
the $\lambda$-ring structure).
\end{enumerate}
The completion of a filtered $\lambda$-ring is again a $\lambda$-ring;
define a complete $\lambda$-ring to be a $\lambda$-ring equal to its
completion. For example, the universal $\lambda$-ring $\Lambda$ is filtered
by the subspaces $F^n\Lambda$ of polynomials vanishing to order $n-1$, and
its completion is the $\lambda$-ring of symmetric power series, whose
underlying ring is the power series ring $\Z\[h_1,h_2,h_3,\dots\]$.

The tensor product of two filtered $\lambda$-rings is again a filtered
$\lambda$-ring, when furnished with the filtration
$$
F^n(R\o S) = \sum_{k=0}^n F^{n-k}R \o F^kS .
$$
If $R$ and $S$ are filtered $\lambda$-rings, denote by $R\ohat S$ the
completion of $R\o S$.

Let $\CR$ be a Karoubian \ring over a field of characteristic zero, and
consider the complete $\lambda$-ring $\Lambda\ohat K_0(\CR)$, where $K_0(\CR)$
has the discrete filtration. This $\lambda$-ring has a natural realization,
as the Grothendieck group of the Karoubian \ring
$$
\hom{\CR} = \prod_{n=0}^\infty [\SS_n,\CR] ,
$$
whose objects are the $\SS$-modules in $\CR$. In this \ring, the sum and
product are given by the same formulas as in the \ring $[\SS,\CR]$ of
bounded $\SS$-modules.

Without assuming the existence of infinite sums in $\CR$, plethysm does not
extend to a monoidal structure on $\hom{\CR}$. However, $\x\circ\y$ is
well-defined in $\[\SS,\CR\]$ under either of the following two hypotheses:
$$
\text{i) $\x$ is bounded, or ii) $\y(0)=0$.}
$$
The first of these situations allows us to construct a $\lambda$-ring
structure on the Grothendieck group of $\hom{\CR}$, by the same method as
for $[\SS,\CR]$, while the second will be needed in the proof of our main
theorem. Introducing the notation $\hom[k]{\CR}$ for the subcategory of
$\hom{\CR}$ consisting of $\SS$-modules $\x$ such that $\x(n)=0$ for $n<k$,
we see that plethysm extends to a symmetric monoidal structure on
$\hom[1]{\CR}$.

Denote the Grothendieck group of the full subcategory
$\hom[1]{\CR}\subset\hom{\CR}$ by $\Check{K}^\SS_0(\CR)$. Since
$\Check{K}^\SS_0(\CR)$ is a (non-unital) $\lambda$-ring, we may define a
bilinear operation
$$
\circ : \Hat{K}^\SS_0(\Proj) \o \Check{K}^\SS_0(\CR) \to \check{K}_\SS(\CR) ,
$$
satisfying \eqref{lambda-ring}. This operation may be extended to a
bilinear operation (which we denote by the same symbol),
$$
\circ : \Hat{K}^\SS_0(\CR) \o \Check{K}^\SS_0(\CR) \to \check{K}_\SS(\CR) ,
$$
using the Peter-Weyl Theorem: to define $x\circ y$, we expand $x$ in a
series $x=\sum_\lambda x_\lambda\*s_\lambda$, and define
$$
x\circ y = \sum_\lambda x_\lambda \* \sigma_\lambda(y) .
$$
The interest of this operation lies in the following lemma, which is a
simple consequence of the definition of $x\circ y$.
\begin{lemma} \label{plethysm}
If $\x$ and $\y$ are objects of $\hom{\CR}$ and $\hom[1]{\CR}$
respectively, $[\x\circ\y]=[\x]\circ[\y]$.
\end{lemma}

If $R$ is a complete $\lambda$-ring, the operation
$$
\Exp(a) = \sum_{n=0}^\infty \sigma_n(a) : R \to 1+F_1R
$$
is an analogue of exponentiation, whose logarithm is given by a formula of
Cadogan \cite{Cadogan}.
\begin{proposition} \label{Cadogan}
On a complete filtered $\lambda$-ring $R$, the operation $\Exp:R\to1+F_1R$
has inverse
$$
\Log(1+a) = \sum_{n=1}^\infty \frac{\mu(n)}{n} \log(1+\psi_n(a)) .
$$
\end{proposition}
\begin{proof}
Expanding $\Log(1+a)$, we obtain
$$
\Log(1+a) = - \sum_{n=1}^\infty \frac{1}{n} \sum_{d|n} \mu(d)
\psi_d(-a)^{n/d} = \sum_{n=1}^\infty \Log_n(a) .
$$
Let $\chi_n$ be the character of the cyclic group $C_n$ equalling $e^{2\pi
i/n}$ on the generator of $C_n$. The characteristic of the $\SS_n$-module
$\Ind^{\SS_n}_{C_n}\chi_n$ equals
$$
\frac{1}{n} \sum_{k=0}^{n-1} e^{2\pi ik/n} p_{(k,n)}^{n/(k,n)}
= \frac{1}{n} \sum_{d|n} \mu(d) p_d^{n/d} ,
$$
while the characteristic of the $\SS_n$-module
$\Ind^{\SS_n}_{C_n}\chi_n\o\eps_n$, where $\eps_n$ is the sign
representation of $\SS_n$, equals
$$
\frac{1}{n} \sum_{d|n} \mu(d) \bigl( (-1)^{d-1}p_d \bigr)^{n/d}
= \frac{(-1)^n}{n} \sum_{d|n} \mu(d) (-p_d)^{n/d} .
$$
It follows that $(-1)^{n-1}\Log_n$ is the operation associated to the
$\SS_n$-module $\Ind^{\SS_n}_{C_n}\chi_n\o\eps_n$, and hence defines a map
from $F_1R$ to $F_nR$.

To prove that $\Log$ is the inverse of $\Exp$, it suffices to check this
for $R=\Lambda$ and $x=h_1$. We must prove that
$$
\Exp\left( \sum_{n=1}^\infty \frac{\mu(n)}{n} \log(1+p_n) \right) = 1+h_1 .
$$
The logarithm of the expression on the left-hand side equals
$$
\exp \Bigl( \sum_{k=1}^\infty \frac{p_k}{k} \Bigr) \circ
\Bigl( \sum_{n=1}^\infty \frac{\mu(n)}{n} \log(1+p_n) \Bigr)
= \sum_{n=1}^\infty \sum_{d|n} \mu(d) \frac{\log(1+p_n)}{n}
= \log(1+p_1) ,
$$
and the formula follows.
\end{proof}

\begin{example} \label{Log(t)}
If $a\in F^1R$ is a line bundle in the complete $\lambda$-ring $R$ (that
is, $\sigma_n(a)=a^n$ for all $n\ge0$), we see that
$$
\Exp(a) = \frac{1}{1-a} .
$$
In particular, this shows that $\Exp(t^n)=(1-t^n)^{-1}$, and that
$$
\Exp(t-t^2) = \frac{\Exp(t)}{\Exp(t^2)} = \frac{1-t^2}{1-t} = 1 + t .
$$
It follows that $\Log(1-t)=t$ and that $\Log(1+t)=t-t^2$.
\end{example}

We now introduce the operations on $\lambda$-rings which will arise in the
calculation of the Serre polynomials of the local systems
$\FF(X,n)\times_{\SS_n}V_\lambda$. We start by considering the case $X=\C$.
\begin{proposition}
$$
\sum_\lambda s_\lambda \o \Serre(\FF(\C,n),V_\lambda)
= \prod_{k=1}^\infty (1+p_k)^{\frac{1}{k}\sum_{d|k}\mu(k/d)\LL^d}
\in \Lambda \ohat \Z[\LL]
$$
\end{proposition}
\begin{proof}
It is proved in Lehrer-Solomon \cite{LS} that
\begin{equation} \label{Lehrer-Solomon}
\sum_{n=0}^\infty \sum_{i=0}^\infty (-x)^i \ch_n(H^i(\FF(\C,n),\C)) =
\prod_{k=1}^\infty (1+x^kp_k)^{\frac{1}{k}\sum_{d|k}\mu(k/d)x^{-d}} ,
\end{equation}
where $H^i(\FF(\C,n),\C)$ is the $\SS_n$-module associated to the de Rham
cohomology of degree $i$. By Poincar\'e duality, we see that
$$
\sum_{n=0}^\infty \sum_{i=0}^\infty (-x)^i \ch_n(H_c^i(\FF(\C,n),\C)) =
\prod_{k=1}^\infty (1+p_k)^{\frac{1}{k}\sum_{d|k}\mu(k/d)x^d} .
$$
But the mixed Hodge structure on the cohomology group $H^i_c(\FF(\C,n),\C)$
is pure of weight $2i$, and indeed
$H^i_c(\FF(\C,n),\C)=H^i_c(\FF(\C,n),\C)^{i,i}$, proving the result.
\end{proof}

Motivated by this proposition, we define operations $\Phi_\lambda$ in a
$\lambda$-ring $R$, parametrized by partitions $\lambda$, by means of the
generating function
\begin{equation} \label{Phi}
\Phi(x) \equiv \sum_\lambda s_\lambda\o\Phi_\lambda(x)
= \prod_{k=1}^\infty (1+p_k)^{\frac{1}{k}\sum_{d|k}\mu(k/d)\psi_d(x)}
\in \Lambda \ohat R .
\end{equation}

\begin{theorem} \label{explicit}
We have the formula $\Phi(x) = \Exp(\Log(1+p_1)x)$. In particular, the
operations $\Phi_\lambda$ are defined on any $\lambda$-ring.
\end{theorem}
\begin{proof}
Applying $\Log$ to the definition of $\Phi(x)$, we obtain
\begin{align*}
\Log(\Phi(x)) &= \sum_{n=1}^\infty \frac{\mu(n)}{n} \psi_n \log(\Phi(x)) \\
&= \sum_{n=1}^\infty \frac{\mu(n)}{n} \psi_n
\sum_{k=1}^\infty \frac{1}{k}\sum_{d|k}\mu(k/d) \log(1+p_k) \psi_d(x) \\
&= \sum_{n,d,e=1}^\infty \frac{\mu(n)\mu(e)}{nde}
\psi_{nd}(x) \log(1+p_{nde}) \\
&= \sum_{e=1}^\infty \frac{\mu(e)}{e} \log(1+p_e) x ,
\end{align*}
by M\"obius inversion. On applying $\Exp$, we obtain the desired formula.
\end{proof}

Using this theorem, we can prove more explicit formulas for $\Phi_n$ and
$\Phi_{1^n}$.
\begin{corollary} \label{braid}
$$
\sum_{n=0}^\infty t^n \Phi_n(y) =
\frac{\sigma_t(y)}{\sigma_{t^2}(y)} \text{ and }
\Phi_{1^n}(y) = \lambda_n(y)
$$
\end{corollary}
\begin{proof}
We obtain $\sum_{n=0}^\infty t^n \Phi_n(y)$ from $\Phi(x)$ by replacing
$p_n$ by $t^n$. By Theorem \ref{explicit}, it follows that
$$
\sum_{n=0}^\infty t^n \Phi_n(y) = \Exp(\Log(1+t)x) = \Exp((t-t^2)x) =
\frac{\sigma_t(x)}{\sigma_{t^2}(x)} ,
$$
since $\Log(1+t)=t-t^2$ by Example \ref{Log(t)}. The proof of the second
formula is similar, except that we replace $p_n$ by $(-t)^n$, and apply the
formula $\Log(1-t)=-t$.
\end{proof}

\section{Representations of finite groups in Karoubian \rings}

Let $(\CR,\otimes,\1)$ be a symmetric monoidal category with coproducts,
denoted $X\oplus Y$. We say that $\CR$ is a \textbf{\ring} (this is our
abbreviation for the usual term \emph{ring category}) if there are natural
isomorphisms
$$
(X\oplus Y)\o Z \cong (X\o Z)\oplus(Y\o Z) \quad\text{and}\quad X\o 0 \cong 0
$$
which describe the distributivity of the tensor product over the sum,
satisfying the coherence axioms of Laplaza \cite{Laplaza}. If $\o$ is the
categorical product, we say that $\CR$ is a Cartesian \ring.

The Grothendieck group $K_0(-)$ is a functor from \rings to commutative
rings. Given an object $X$ of a \ring $\CR$, denote by $[X]$ its
isomorphism class; then $K_0(\CR)$ is generated as an abelian group by the
isomorphism classes of objects, with the relation
$$
[X] + [Y] = [X\oplus Y] .
$$
The product on $K_0(\CR)$ is given by the formula $[X]\*[Y]=[X\o
Y]$. (Here, we suppose that the isomorphism classes of objects of $\CR$
form a set; this hypothesis will always be fulfilled in this paper.)

If $\CR$ and $\CS$ are two \rings, $\CR\o\CS$ is a \ring whose objects
are formal sums of tensor products $X\o Y$, where $X$ and $Y$ are objects
of $\CR$ and $\CS$ respectively; note that $K_0(\CR\o\CS)\cong K_0(\CR)\o
K_0(\CS)$.

Recall that an additive category over a commutative ring $R$ is a category
$\CR$ such that the set of morphisms $\CR(X,Y)$ is a $R$-module for all
objects $X$ and $Y$, the composition maps $\CR(Y,Z)\o_K\CR(X,Y)\to\CR(X,Z)$
are $R$-linear, and every finite set of objects has a direct sum. A
\textbf{Karoubian category} over a ring $R$ is an additive category over
$R$ such that every idempotent has an image, denoted $\im(p)$. (Karoubian
categories are also sometimes known as pseudo-abelian categories.)
\begin{definition}
A \textbf{Karoubian \ring} $\CR$ is a \ring which is a Karoubian category,
and whose sum $X\oplus Y$ is the direct sum.
\end{definition}

An example of a Karoubian \ring is the category $\Proj$ of finitely
generated projective $R$-modules.

If $\CR$ is a Karoubian \ring and $G$ is a group, let $[G,\CR]$ be the
Karoubian \ring of $G$-modules in $\CR$, that is, functors from $G$ to
$\CR$. If $X$ and $Y$ are objects of $[G,\CR]$, the $R$-module of morphisms
$\CR(X,Y)$ carries a natural $R[G]$-module structure, given by the formula
$f^g=g^{-1}\*f\*g$.

There is a natural bifunctor $V\boxtimes X$, the external tensor product,
from $[G,\Proj]\times[G,\CR]$ to $[G,\CR]$, characterized by the identity
of $R[G]$-modules
$$
\CR(V\boxtimes X,Y) \cong V \o \CR(X,Y) .
$$
For the finitely generated free module $R[G]^n$, we have
$$
R[G]^n\boxtimes X = \bigoplus_{g\in G} X^{\oplus n} .
$$
For general $V$, we realize $V$ as the image of an idempotent $p$ in a free
module $R[G]^n$, and define $V\boxtimes X$ to be the image of the
corresponding idempotent in $R[G]^n\boxtimes X$. Using the external tensor
product, we may embed $[G,\Proj]$ into $[G,\CR]$ by the functor $V\mapsto
V\boxtimes\1$.

If $G$ is a group whose order is invertible in $R$, the functor $(-)^G$ of
$G$-invariants from $[G,\CR]$ to $\CR$ is defined by taking the image of
the idempotent automorphism of $\CR$
$$
p = \frac{1}{|G|} \sum_{g\in G} g .
$$
From now on, we restrict attention to groups satisfying this condition.

If $H$ is a subgroup of $G$, the induction functor
$\Ind^G_H:[H,\CR]\to[G,\CR]$ is defined by the formula
$$
\Ind^G_H X = (R[G]\boxtimes X)^H .
$$
Here, we use the $G\times H$-module structure of $R[G]$, where $G$ acts on
the left and $H$ acts on the right.

The following is a generalization of the Peter-Weyl theorem to Karoubian
categories.
\begin{theorem}[Peter-Weyl] \label{Peter-Weyl}
If $\CR$ is a Karoubian \ring over a commutative ring $R$ and $G$ is a
group whose order is invertible in $R$, the composition
$$
[G,\Proj] \o \CR \hookrightarrow [G,\Proj] \o [G,\CR]
\xrightarrow{\boxtimes} [G,\CR]
$$
is an equivalence of categories.
\end{theorem}
\begin{proof}
Since the order of $G$ is invertible in $R$, the group algebra $R[G]$ is
semi-simple, and may be written
$$
R[G] \cong \bigoplus_a \End(V_a) \cong \bigoplus_a V_a \o V_a^* ,
$$
where we sum over the isomorphism classes of irreducible representations
$\{V_a\}$ of $G$. This permits us to rewrite the induction functor as
$$
\Ind_H^GX = (R[G]\o X)^H \cong \bigoplus_a V_a \o (V_a^*\boxtimes X)^H .
$$
Taking $H=G$, and recalling that $\Ind_G^G$ is equivalent to the identity
functor, we obtain the desired equivalence between $[G,\Proj]\o\CR$ and
$[G,\CR]$.
\end{proof}

\section{$\SS$-modules in Karoubian \rings}

Let $\SS$ be the category of permutations $\coprod_{n=0}^\infty\SS_n$ and
let $\CR$ be a \ring. A bounded $\SS$-module in $\CR$ is an object $\x$ of
$$
[\SS,\CR] = \bigoplus_{n=0}^\infty \, [\SS_n,\CR] ,
$$
in other words, a sequence $\{\x(n)\mid n\ge0\}$ of $\SS_n$-modules in
$\CR$ such that $\x(n)=0$ for $n\gg0$. Let $\1_n$ denote the $\SS$-module
such that $\1_n(n)$ is the trivial $\SS_n$-module, while $\1_n(k)=0$ for
$k\ne n$.

The category $[\SS,\CR]$ is itself a \ring:
\begin{enumerate}
\item the sum of two $\SS$-modules is $(\x\oplus\y)(n)=\x(n)\oplus\y(n)$;
\item the product of two $\SS$-modules is defined using induction:
$$
(\x\o\y)(n) = \bigoplus_{j+k=n} \Ind_{\SS_j\times\SS_k}^{\SS_n} \x\o\y ;
$$
\item the unit of the product is $\1_0$.
\end{enumerate}
Denote the Grothendieck group of the \ring $[\SS,\CR]$ by $K_0^\SS(\CR)$.

There is another monoidal structure $\x\circ\y$ on $[\SS,\CR]$, called
plethysm. If $\lambda$ is a partition of $n$, let $\SS_\lambda =
\SS_{\lambda_1}\times\dots\times\SS_{\lambda_{\ell(\lambda)}}\subset\SS_n$,
and let $N(\SS_\lambda)$ be the normalizer of $\SS_\lambda$ in $\SS_n$. The
quotient $W(\SS_\lambda)=N(\SS_\lambda)/\SS_\lambda$ may be identified with
$$
\{ \sigma\in\SS_{\ell(\lambda)} \mid \text{$\lambda_{\sigma(i)}=\lambda_i$
for $1\le i\le\ell(\lambda)$} \} \subset \SS_{\ell(\lambda)} .
$$
Given bounded $\SS$-modules $\x$ and $\y$, we obtain an action of
$N(\SS_\lambda)$ on the tensor product
$$
\x(\ell(\lambda)) \o \bigotimes_{1\le i\le\ell(\lambda)} \y(\lambda_i) .
$$
\textbf{Plethysm} is the monoidal structure (not symmetric) defined by the
formula
$$
(\x\circ\y)(n) = \bigoplus_{\lambda\vdash n} \bigoplus_{k=0}^\infty
\Ind^{\SS_n}_{N(\SS_\lambda)}
\biggl( \x(\ell(\lambda)+k) \o \bigotimes_{1\le i\le\ell(\lambda)}
\y(\lambda_i) \o \y(0)^{\o k} \biggr)^{\SS_k} ,
$$
and with unit $\1_1$.

\begin{lemma}
Let $\CR$ be a Karoubian \ring over a field of characteristic zero. The
Grothendieck group $K_0^\SS(\CR)$ is a pre-$\lambda$-ring, with
$\sigma$-operations characterized by the formula
$$
\sigma_n([\x]) = \bigl[\1_n\circ\x\bigr] ,
$$
where $\x$ is a bounded $\SS$-module.
\end{lemma}
\begin{proof}
We must prove that for bounded $\SS$-modules $\x$ and $\y$,
\begin{equation} \label{pre}
\sigma_n([\x]+[\y]) = \sum_{i=0}^n \sigma_i([\x]) \* \sigma_{n-i}([\y]) .
\end{equation}
Observe that $\1_n\circ(\x\oplus\y)$ equals
$$
\bigoplus_{i=0}^n \bigoplus_{\lambda\vdash i} \bigoplus_{\mu\vdash n-i}
\bigoplus_{j,k=0}^\infty \Ind^{\SS_n}_{N(\SS_\lambda)\times N(\SS_\mu)}
\biggl( \bigotimes_{1\le i\le\ell(\lambda)} \x(\lambda_i) \o \x(0)^{\o j}
\o \bigotimes_{1\le i\le\ell(\mu)} \y(\mu_i) \o \y(0)^{\o k}
\biggr)^{\SS_j\times\SS_k}
$$
Since
$$
\Ind^{\SS_n}_{N(\SS_\lambda)\times N(\SS_\mu)}V\o W
= \Ind^{\SS_n}_{\SS_i\times\SS_{n-i}} \Bigl(
\Ind^{\SS_i}_{N(\SS_\lambda)} V \o \Ind^{\SS_{n-i}}_{N(\SS_\mu)} W \Bigr) ,
$$
it follows that
$$
\1_n\circ(\x\oplus\y) \cong \bigoplus_{i=0}^n
(\1_i\circ\x)\o(\1_{n-i}\circ\y) ,
$$
proving \eqref{pre} for elements of $K_0^\SS(\CR)$ of the form $[\x]$ and
$[\y]$. The definition of the sigma operations on virtual elements
$[\x_0]-[\x_1]$ is now forced by \eqref{invert}:
$$
\sigma_n([\x_0]-[\x_1]) = \sum_{k=0}^\infty 
\sum_{\substack{j_\ell>0 \\ i+j_1+\dots+j_\ell=n}}
(-1)^k \sigma_i([\x_0]) \sigma_{j_1}([\x_1]) \dots \sigma_{j_k}([\x_1]) .
\qed$$
\def\qed{}
\end{proof}

\begin{lemma} \label{atiyah}
There is an isomorphism of $\lambda$-rings $K_0^\SS(\Proj)\cong\Lambda$.
\end{lemma}
\begin{proof}
The pre-$\lambda$-ring $K_0^\SS(\Proj)$ is the sum of abelian groups
$
K_0^\SS(\CR) = \bigoplus_{n=0}^\infty R(\SS_n) ,
$
where $R(\SS_n)=K_0([\SS_n,\Proj)$ is the abelian group underlying the
virtual representation ring of $\SS_n$. The identification of
$K_0^\SS(\Proj)$ with $\Lambda$ is via the Frobenius characteristic
$\ch:R(\SS)\to\Lambda$, which sends the irreducible representation
$V_\lambda$ associated to the partition $\lambda$ to the Schur function
$s_\lambda$. The Frobenius characteristic is given by the explicit formula
$$
\ch_n(V) = \frac{1}{n!} \sum_{\sigma\in\SS_n} \Tr_V(\sigma) p_\sigma ,
$$
where $p_\sigma$ is the monomial in the power sums obtained by taking one
factor $p_k$ for each cycle of $\sigma$ of length $k$. For the proof that
$\ch(\dots)$ is a map of $\lambda$-rings, see Knutson \cite{Knutson} or
Appendix~A of Macdonald \cite{Macdonald}.
\end{proof}

Using these lemmas and the Peter-Weyl Theorem, we will show that
$K_0^\SS(\CR)$ is a $\lambda$-ring for any Karoubian \ring over a field of
characteristic zero. First, we prove some simple lemmas which are of
interest in their own right.

Plethysm is distributive on the left with respect to sum.
\begin{lemma} \label{additive}
$
(\x_1\oplus\x_2)\circ\y \cong (\x_1\circ\y) \oplus (\x_2\circ\y)
$
\end{lemma}
\begin{proof}
Clear.
\end{proof}

It is also distributive on the left with respect to product.
\begin{lemma} \label{multiplicative}
$
(\x_1\o\x_2)\circ\y \cong (\x_1\circ\y) \o (\x_2\circ\y)
$
\end{lemma}
\begin{proof}
By Lemma \ref{additive}, it suffices to check this formula when
$\x_1(j)=X_1$, $\x_2(k)=X_2$, $\x_1(i)=0$, $i\ne j$ and $\x_1(i)=0$, $i\ne
k$. We have
\begin{multline*}
\bigl((\x_1\o\x_2) \circ \y\bigr)(n) \\
= \bigoplus_{q=0}^n
\bigoplus_{\substack{\lambda\vdash n \\ \ell(\lambda)+q=j+k}}
\Ind^{\SS_n}_{N(\SS_\lambda)} \biggl( \Ind^{\SS_{j+k}}_{\SS_j\times\SS_k}
\bigl( X_1\o X_2 \bigr) \o \bigotimes_{1\le i\le\ell(\lambda)}
\y(\lambda_i) \o \y(0)^{\o q} \biggr)^{\SS_q} .
\end{multline*}
But we have
\begin{multline*}
\bigoplus_{\substack{\lambda\vdash n \\ \ell(\lambda)+q=j+k}} \biggl(
\Ind^{\SS_{j+k}}_{\SS_j\times\SS_k} \bigl( X_1\o X_2 \bigr) \o
\bigotimes_{1\le i\le\ell(\lambda)} \y(\lambda_i) \o \y(0)^{\o q}
\biggr)^{\SS_q} \cong \bigoplus_{p=0}^q \bigoplus_{i=0}^n \\
\bigoplus_{\substack{\mu\vdash i \\ \ell(\mu)+p=j}}
\biggl( X_1 \o \bigotimes_{1\le i\le\ell(\mu)} \y(\lambda_i) \o
\y(0)^{\o p} \biggr)^{\SS_q}
\o \bigoplus_{\substack{\lambda\vdash n-i \\ \ell(\lambda)+q-p=k}}
\biggl( X_2 \o \bigotimes_{1\le i\le\ell(\lambda)} \y(\lambda_i) \o
\y(0)^{\o q} \biggr)^{\SS_q} ,
\end{multline*}
from which the lemma follows.
\end{proof}

\begin{lemma} \label{associative}
If $\V$ is a bounded $\SS$-module in $\Proj$ and $\x$ is a bounded
$\SS$-module in $\CR$,
$$
\ch(\V) \circ [\x] = [\V\circ\x] .
$$
\end{lemma}
\begin{proof}
By Lemma \ref{additive}, we may assume that $\V$ is an irreducible
$\SS_n$-module $V_\lambda$. It remains to show that
$\sigma_\lambda([\x])=[V_\lambda\circ\x]$ for all partitions $\lambda$.

By Lemma \ref{multiplicative}, we see that for any partition $\mu$ with
$\ell=\ell(\mu)$, we have
$$
\bigl( \1_{\mu_1} \o \dots \o \1_{\mu_\ell} \bigr) \circ \x
\cong \bigl( \1_{\mu_1} \o \x \bigr) \o \dots \o
\bigl( \1_{\mu_\ell} \circ \x \bigr) .
$$
Taking the class in $K_0^\SS(\CR)$ of both sides, we see that
$$
\bigl[ \bigl( \1_{\mu_1} \o \dots \o \1_{\mu_\ell} \bigr) \circ \x \bigr] =
\sigma_{\mu_1}([\x]) \dots \sigma_{\mu_\ell}([\x]) .
$$
The irreducible representation $V_\lambda$ is a linear combination of
representations $\1_{\mu_1}\o\dots\o\1_{\mu_\ell}$ with integral
coefficients, and by Lemma \ref{atiyah}, the Schur function $s_\lambda$ is
a linear combination of symmetric functions $h_{\mu_1}\o\dots\o
h_{\mu_\ell}$ with the same coefficients; the proof is completed by
application of Lemma \ref{additive}.
\end{proof}

\begin{theorem}
The Grothendieck group $K_0^\SS(\CR)$ of a Karoubian \ring $\CR$ over a field
of characteristic zero is a $\lambda$-ring.
\end{theorem}
\begin{proof}
If $f=\ch(\V)$ and $g=\ch(\w)$, where $\V$ and $\w$ are bounded
$\SS$-modules in $\Proj$, and $x=[\x]$, where $\x$ is a bounded
$\SS$-module in $\CR$, it follows from Lemma \ref{associative} that
$$
f \circ \bigl( g \circ x \bigr) = \ch\bigl( \V \circ (\w\circ\x) \bigr)
= \ch\bigl( (\V\circ\w) \circ \x \bigr) = \ch(\V\circ\w) \circ x .
$$
Since $\ch$ is a morphism of $\lambda$-rings, we see that
$\ch(\V\circ\w)=f\circ g$, and from which we obtain the formula
\eqref{lambda-ring} characterizing $\lambda$-rings in this case:
$$
f\circ(g\circ x) = (f\circ g) \circ x .
$$
It only remains to extend \eqref{lambda-ring} to virtual elements
$g=\ch(\w_0)-\ch(\w_1)$ and $x=[\x_0]-[\x_1]$. Both sides of
\eqref{lambda-ring} are polynomial functions of $g\in\Lambda$ and $x\in
K_0^\SS(\CR)$ and hence must coincide, since they are equal on a cone with
non-empty interior.
\end{proof}

It follows that the Grothendieck group $K_0(\CR)$ is a $\lambda$-ring,
namely the sub-$\lambda$-ring of $K_0^\SS(\CR)$ consisting of virtual
objects such that $X(n)=0$, $n>0$. The Peter-Weyl Theorem now has the
following consequence.
\begin{theorem}
If $\CR$ is a Karoubian \ring over a field of characteristic zero, there is
an isomorphism $K_0^\SS(\CR)\cong\Lambda\o K_0(\CR)$ of $\lambda$-rings.
\end{theorem}
\begin{proof}
The Peter-Weyl Theorem gives isomorphisms of rings
$$
\Lambda \o K_0(\CR) \xleftarrow{\ch\o1} K_0^\SS(\Proj) \o K_0(\CR)
\xrightarrow{\boxtimes} K_0^\SS(\CR) .
$$
The first of these arrows is an isomorphism of $\lambda$-rings by Lemma
\ref{atiyah}. As rings, both $K_0^\SS(\Proj)\o K_0(\CR)$ and $K_0^\SS(\CR)$
are generated by $K_0(\CR)$ and $[\1_n]$, $n\ge1$, and $\boxtimes$ respects
the $\sigma$-operations of these elements, proving that it is a map of
$\lambda$-rings.
\end{proof}

\section{The main result}

If $\CR$ is a Karoubian \ring, denote by $\CR[\N]$ the Karoubian ring of
bounded sequences
$$
(A^0,A^1,A^2,\dots\mid\text{$A^n=0$ for $n\gg0$}) .
$$
The sum on $\CR[\N]$ is defined by $(A\oplus\y)^n=A^n\oplus\y^n$, while the
product is defined by
$$
(A\o B)^n = \bigoplus_{i+j=n} A^i \o B^j .
$$
\begin{definition}
A \textbf{K\"unneth functor} with values in the Karoubian \ring $\CR$ is a
\ring functor $\GS$ from the Cartesian {\ring} $\Var$ of quasi-projective
varieties and open embeddings to $\CR[\N]$.
\end{definition}

In other words, a functor $\GS:\Var\to\CR[\N]$ is a K\"unneth functor if
there are natural isomorphisms
\begin{gather*}
\textstyle
\GS^i(X\coprod Y)\cong\GS^i(X)\oplus\GS^i(Y) , \\
\GS^n(X\times Y) \cong \bigoplus_{n=i+j} \GS^i(X)\o\GS^j(Y) .
\end{gather*}

If $\GS=\{\GS^n\}$ is a K\"unneth functor, denote by $\Serre(X)$ the
associated Euler characteristic
$$
\Serre(X) = \sum_{n=0}^\infty (-1)^n [\GS^n(X)]
$$
in the Grothendieck group $K_0(\CR)$.

\begin{definition}
A \textbf{Serre functor} with values in $\CR$ is a K\"unneth functor $\GS$
such that for any closed sub-variety $Z$ of $X$,
$$
\Serre(X) = \Serre(X\setminus Z) + \Serre(Z) .
$$
\end{definition}

If $\GS$ is a Serre functor and $X=X^0\subset X^1\subset X^2\subset\dots$
is a filtered quasi-projective variety such that $X^n=\emptyset$ for
$n\gg0$, we have
\begin{equation} \label{gr}
\Serre(\gr X) \equiv \sum_n \Serre(\gr^nX) = \Serre(X) .
\end{equation}

Here are two examples of Serre functors:
\begin{enumerate}
\item The category of mixed Hodge structures over $\C$ is a \ring, whose
Grothendieck group may be identified with the polynomial ring $\Z[u,v]$ by
means of the Serre polynomial \eqref{Serre}. The functor $\GS^n(X)$ which
takes a quasi-projective variety $X$ to the mixed Hodge structure
$(H^n_c(X,\C),F,W)$ over $\C$ is a Serre functor. The associated
characteristic $\Serre(X)$ may be identified with the Serre polynomial.
\item Gillet and Soul\'e \cite{GS} have constructed a functor to the
homotopy category of chain complexes of (pure effective rational) Chow
motives; let $\GS^n(X)$ be the $n$th cohomology of this complex.
\end{enumerate}

If $\CR$ is a \ring, let $\T:\CR\to\[\SS,\CR\]$ be the \ring functor with
$\T(X,n)=X^n$. (More precisely, $\T(X,n)$ is defined by induction:
$\T(X,0)=\1$, and $\T(X,n)=\T(X,n-1)\o X$.) The following result is a
generalization of Macdonald's formula \cite{Macdonald-symmetric} for the
Poincar\'e polynomial of the symmetric power $S^nX=X^n/\SS_n$.
\begin{proposition} \label{Macdonald}
If $X$ is a quasi-projective variety,
$$
\Serre(\T(X)) = \Exp\bigl(p_1\Serre(X)\bigr) \in \Hat{K}^\SS_0(\CR) .
$$
Here, $\Serre(\T(X))$ denotes the class $n\mapsto\Serre(\T(X,n))$ in the
Grothendieck group $\Hat{K}^\SS_0(\CR)$.
\end{proposition}
\begin{proof}
Since $\GS$ is a \ring-functor, $\GS\*\T=\T\*\GS$. By the Peter-Weyl
Theorem,
$$
\GS(\T(X,n)) = \T[\GS(X)](n)
= \bigoplus_{\lambda\vdash n} V_\lambda \boxtimes \bigl( V_\lambda^*\o
\GS(X)^{\o n} \bigr)^{\SS_n} .
$$
Descending to the Grothendieck group, we see that
$$
\Serre(\T(X,n)) = \bigoplus_{\lambda\vdash n} s_\lambda \o
\sigma_\lambda(\Serre(X)) \in \Lambda_n \o K_0(\CR) \subset K_0^\SS(\CR) .
$$
Summing over $n\ge0$, and applying Cauchy's formula \eqref{Cauchy}, we see
that
$$
\Serre(\T(X)) = \exp \Bigl( \sum_{k=1}^\infty
\frac{p_k\o\psi_k\Serre(X)}{k} \Bigr) \in \Hat{K}^\SS_0(\CR) .
$$
The proposition now follows by the definition of $\Exp(\dots)$.
\end{proof}

Consider the following decreasing filtration on the $\SS$-module $\T(X)$,
where $X$ is a quasi-projective variety:
$$
\T^i(X)(n) = \{ (z_1,\dots,z_n) \in X^n \mid \text{$\{z_1,\dots,z_n\}$
has cardinality at most $n-i$} \} .
$$
Let $\gr^i\T(X)=\T^i(X)\setminus\T^{i+1}(X)$ be the associated graded
$\SS$-module.
\begin{lemma} \label{gr-T}
Let $\ZZ$ be the object of $\hom[1]{\Var}$
$$
\ZZ(n) = \begin{cases} \point , & n>0 , \\ \emptyset , & n=0 .
\end{cases}$$
Then $\gr\T(X)=\FF(X)\circ\ZZ$; in particular, $\gr^0\T(X)=\FF(X)$.
\end{lemma}
\begin{proof}
This lemma reflects the fact that an element of $\gr^i\T(X,n)$ determines,
and is determined by, a partition of the set $\{1,\dots,n\}$ into $n-i$
disjoint subsets, together a point in $\FF(X,n-i)$.
\end{proof}

We now arrive at the main theorem of this paper.
\begin{theorem} \label{MAIN}
Let $X$ be a quasi-projective variety over $\C$. If $\GS$ is a Serre
functor and $V_\lambda$ is an irreducible representation of $\SS_n$,
$$
\Serre(\FF(X,n),V_\lambda) = \Phi_\lambda(\Serre(X)) .
$$
\end{theorem}
\begin{proof}
If $\GS$ is Serre functor, \eqref{gr} and Lemma \ref{gr-T} show that
$$
\Serre(\T(X)) = \Serre(\gr\T(X)) = \Serre(\FF(X)) \circ \Serre(\ZZ) .
$$
To calculate $\Serre(\FF(X))$, we invert the operation $-\circ\Serre(\ZZ)$
on $\Hat{K}^\SS_0(\CR)$. Indeed, $\Serre(\ZZ)=\Exp(p_1)-1$ and by Lemma
\ref{Cadogan},
\begin{align*}
\Serre(\FF(X)) &= \Serre(\FF(X)) \circ \bigl( \Exp(p_1) - 1 \bigr)
\circ \Bigl( \sum_{n=1}^\infty \frac{\mu(n)}{n} \log(1+p_n) \Bigr) \\
&= \Serre(\T(X)) \circ \Bigl( \sum_{n=1}^\infty \frac{\mu(n)}{n}
\log(1+p_n) \Bigr) .
\end{align*}
By Proposition \ref{Macdonald}, this equals
\begin{multline*}
\exp \Bigl( \sum_{k=1}^\infty \frac{p_k\*\psi_k\Serre(X)}{k} \Bigr) \circ
\Bigl( \sum_{\ell=1}^\infty \frac{\mu(\ell)}{\ell} \log(1+p_\ell) \Bigr)
= \exp \Bigl( \sum_{k=1}^\infty \sum_{\ell=1}^\infty \frac{\mu(\ell)}{k\ell}
\log(1+p_{k\ell}) \psi_k\Serre(X)\Bigr) \\
= \exp \Bigl( \sum_{n=1}^\infty \sum_{d|n} \frac{\mu(n/d)}{n}
\log(1+p_n) \psi_d\Serre(X) \Bigr) ,
\end{multline*}
from which the theorem follows by extracting the coefficient of the Schur
function $s_\lambda$ on both sides.
\end{proof}

The concise formulation
$$
\Serre(\FF(X)) = \Exp(\Log(1+p_1)\Serre(X))
$$
of this result makes the resemblance with the formula
$\Serre(\T(X))=\Exp(p_1\Serre(X))$ clearer.

In the special cases $\lambda=(n)$ or $\lambda=(1^n)$, when $\Phi_\lambda$
is given by the explicit formula of Corollary \ref{braid}, we obtain the
following corollary.
\begin{corollary}
If $\Serre(X)=\sum_{p,q}h_{pq}u^pv^q$ is the Serre polynomial of $X$, then
\begin{gather*}
\sum_{n=0}^\infty t^n \Serre(\FF(X,n)/\SS_n) =
\frac{\sigma_t(X)}{\sigma_{t^2}(X)} = \prod_{p,q=0}^\infty
\left( \frac{1-t^2u^pv^q}{1-tu^pv^q} \right)^{h_{pq}} , \\
\sum_{n=0}^\infty t^n \Serre(\FF(X,n),\eps) = \sigma_{-t}(X)^{-1} =
\prod_{p,q=0}^\infty (1+tu^pv^q)^{h_{pq}} .
\end{gather*}
\end{corollary}

For example, if $X=\C$, $\FF(\C,n)/SS_n$ is the classifying space
$K(B_n,1)$ of the braid group $B_n$ on $n$ strands. Our formula becomes in
this case
$$
\sum_{n=0}^\infty t^n \Serre(\FF(\C,n)/\SS_n) = \frac{1-t^2uv}{1-tuv}
= 1 + t\LL + t^2(\LL^2-\LL) + t^3(\LL^3-\LL^2) + \dots ,
$$
reflecting the isomorphism of rational cohomology groups
$H^\bull(B_n,\Q)\cong H^\bull(\mathbb{G}_m,\Q)$ as mixed Hodge structures.

\section{The Fulton-MacPherson compactification}

Fulton and MacPherson \cite{FM} have introduced a sequence of functors
$X\mapsto X[n]$ from $\Var$ to $[\SS_n,\Var]$, with the following
properties.
\begin{enumerate}
\item If $X$ is projective, then so is $X[n]$.
\item There is natural transformation of functors $\FF(X,n)\hookrightarrow
X[n]$, which is an embedding.
\item The complement $X[n]\setminus\FF(X)$ is a divisor with normal
crossings.
\end{enumerate}
In this section, we calculate the equivariant Serre polynomial
$\Serre(X[n])$. Denote by $\FM(X)$ the functor $X\mapsto(n\mapsto X[n])$
from $\Var$ to $\[\SS,\Var\]$,

\subsection{Trees and $\SS$-modules}
Let $\Gamma(n)$, $n\ge2$, be the set of isomorphism classes of labelled
rooted trees with $n$ leaves, such that each vertex has at least two
branches. It is easily seen that $\Gamma(n)$
is finite: in fact, the generating function
\begin{equation} \label{enumerate}
x + \sum_{n=2}^\infty \frac{x^n |\Gamma(n)|}{n!}
\end{equation}
is the inverse under composition of $x-x^2-x^3-x^4-\dots$.

Given a tree $T\in\Gamma(n)$, denote by $\VERT(T)$ the set of vertices of
$T$; given a vertex $v\in\VERT(T)$, denote by $n(v)$ the valence of $v$
(its number of branches). Given a tree $T\in\Gamma(n)$ and an $\SS$-module
$\V$ in the \ring $\CR$, let $\V(T)$ be the object
\begin{equation} \label{V(T)}
\V(T) = \bigotimes_{v\in\VERT(T)} \V(n(v)) ,
\end{equation}
and let $\TT\V(n)$ be the $\SS_n$-module
$$
\TT\V(n) =  \bigoplus_{T\in\Gamma(n)} \V(T) .
$$
Thus, $\TT$ is a functor from $\hom[2]{\CR}$ to itself. (Recall that
$\hom[2]{\CR}$ is the full subcategory of $\SS$-modules such that
$X(0)=X(1)=0$.)

A proof of the following formula for $\CR=\Proj$ may be found in
\cite{modular}; however, the same proof works in general. Observe that this
theorem may be used to prove \eqref{enumerate}.
\begin{theorem} \label{revert}
The elements
$$
f = h_1 - \sum_{n=2}^\infty [\V] \quad\text{and}\quad
g = h_1 + \sum_{n=2}^\infty [\TT\V]
$$
of $\Check{K}^\SS_0(\CR)$ satisfy the formula $f\circ g = g\circ f=h_1$.
\end{theorem}

\subsection{The varieties $\protect\Po_k(n)$}
The algebraic groups $\C^k$ and $\G_m$ act on the affine space $\C^k$ by
translation and dilatation respectively; by functoriality, these actions
extend to $\FF(\C^k,n)$. Denote by $G_k=\C^k\rtimes\G_m$ the semidirect
product of these groups, and by $\Po_k(n)$, $n>1$, the quotient of the
configuration space $\FF(\C^k,n)$ by the free $G_k$-action. This action is
$\SS_n$-equivariant, and $\Po_k(n)$ is a smooth $\SS_n$-variety of
dimension $nk-k-1$. For example, $\Po_k(2)$ is naturally isomorphic to the
projective space $\CP^{k-1}$, with trivial $\SS_2$-action.

\begin{proposition}
$$
\Serre\bigl( \Po_k(n),\SS_n \bigr) =
\frac{\Serre(\FF(\C^k,n),\SS_n)}{\Serre(\C^k)\Serre(\G_m)} =
\frac{\Serre(\FF(\C^k,n),\SS_n)}{\LL^k(\LL-1)}
$$
\end{proposition}
\begin{proof}
We start with a lemma.
\begin{lemma}
Let $G$ be an algebraic group and $P$ be a $G$-torsor with base $X=P/G$.
If the projection $P\to X$ is locally trivial in the Zariski topology,
$\Serre(P) = \Serre(G) \Serre(X)$.
\end{lemma}
\begin{proof}
We stratify $X$ by locally closed subvarieties $X_i$ of codimension $i$
over which the torsor $P$ is trivial. The strata are chosen inductively:
$X_{-1}$ is empty, while $X_i$ is a Zariski-open subset of $X\setminus
X_{i-1}$ over which $P$ is trivial. The formula follows, since
$$
\Serre(P) = \sum_i \Serre(P_i) = \sum_i \Serre(G) \Serre(X_i) .
\qed$$
\def\qed{}
\end{proof}

The action of $\C^k$ on $\FF(\C^k,n)$ is not just locally, but globally,
trivial: a global section is given by $(z_1,\dots,z_n) \mapsto
(z_1-\zbar,\dots,z_n-\zbar)$, where $\zbar = \frac{1}{n} \sum_{i=1}^n z_i$.
On the other hand, any free action of $\G_m$ on a variety is locally
trivial in the Zariski topology: free actions with quotient $X$ are
classified by $H^1(X_{\text{fl}},\G_m)$, locally trivial free actions with
quotient $X$ are classified by $H^1(X,\G_m)$, and these two groups are
isomorphic by Hilbert's Theorem 90 (see Proposition XI.5.1 of Grothendieck
\cite{SGA1}).
\end{proof}

\subsection{Stratification of $\FM(X)$} 
The $\SS$-variety $\Po_k$ has a natural compactification to a smooth
projective $\SS$-variety $\P_k$, which has a natural stratification. The
strata are labelled by trees $T\in\Gamma(n)$, and the stratum associated to
$T$ is isomorphic to $\Po_k(T)$, in the notation of \eqref{V(T)}. It
follows from Theorem \ref{revert} that $\Serre(\P_k)$ is the inverse of
\begin{align*}
h_1 - \Serre\bigl( \Po_k \bigr) &= p_1 - \frac{\prod_{n=1}^\infty
(1+p_n)^{\frac{1}{n}\sum_{d|n}\mu(n/d)\LL^{kd}} - 1 - \LL^kp_1}{\LL^k(\LL-1)} \\
& = \frac{\LL^{k+1}p_1 + 1 - \prod_{n=1}^\infty
(1+p_n)^{\frac{1}{n}\sum_{d|n}\mu(n/d)\LL^{kd}}}{\LL^k(\LL-1)} .
\end{align*}
under plethysm.

The importance of the spaces $\P^k(n)$ comes from the following result of
Fulton and MacPherson.
\begin{proposition}
The $\SS$-module $\FM(X)$ has a filtration such that
$$
\gr\FM(X) \cong \FF(X)\circ\P_k .
$$
\end{proposition}

Since $X[n]$ is a projective $Q$-variety (it has singularities which are
quotients of affine space by a finite group), $\Serre(\FM(X))(n)$ equals
the $\SS_n$-equivariant Hodge polynomial of $X[n]$. The above proposition
shows that $\Serre(\FM(X))=\Serre(\FF(X))\circ\Serre(\P_k)$, and leads to a
practical algorithm for the calculation of the $\SS_n$-equivariant Hodge
numbers of $X[n]$.

On forgetting the action of the symmetric groups $\SS_n$, we recover the
formula of Fulton and Macpherson for the Poincar\'e polynomials of
$\FM(X,n)$, in a form stated by Manin \cite{Manin:1}. On replacing $h_n$
by $x^n/n!$, we obtain
$$
1 + \sum_{n=1}^\infty x^n \Serre(X[n]) = (1+x)^{\Serre(X)} \circ
\biggl( \frac{\LL^{k+1}x+1-(1+x)^{\LL^k}}{\LL^k(\LL-1)} \biggr)^{-1}
$$
In this formula, we may take the limit $\LL\to1$ using L'H\^opital's rule,
obtaining a formula for the Euler characteristic of $\FM(X,n)$:
$$
1 + \sum_{n=1}^\infty x^n \chi(X[n])
= (1+x)^{\chi(X)} \circ \bigl( (k+1)x - k(1+x)\log(1+x) \bigr)^{-1} .
$$

The one dimensional case has special interest, since the spaces $\Po_1(n)$
and $\P_1(n)$ are naturally isomorphic to the moduli spaces $\CM_{0,n+1}$
and $\Mbar_{0,n+1}$; this isomorphism comes about because the translations
and dilatations in one dimension generate the isotropy group of the point
$\infty\in\CP^1$ with respect to the action of the group $\PSL(2,\C)$. This
identification means that the action of $\SS_n$ on these spaces is the
restriction of an action of $\SS_{n+1}$. We have calculated the
$\SS_{n+1}$-equivariant Serre polynomials of these spaces in
\cite{gravity}; in a sequel to this paper, we calculate the
$\SS_n$-equivariant Serre polynomial of $\Mbar_{1,n}$.


\end{document}